\documentstyle[12pt,amsfonts]{article}

\newcommand{\theoremName}{Theorem}
\newcommand{\lemmaName}{Lemma}
\newcommand{\corollaryName}{Corollary}
\newcommand{\statementName}{Statement}
\newcommand{\remarkName}{Remark}
\newcommand{\theremark}{\arabic{remark}}
\newcommand{\exampleName}{Example}
\newcommand{\theexample}{\arabic{example}}
\newcommand{\definitionName}{Definition}
\newcommand{\thedefinition}{\arabic{definition}}
\newcommand{\problemName}{Problem}
\newcommand{\theproblem}{\arabic{problem}}

\newtheorem {theorem}{\theoremName}
\newtheorem {lemma}{\lemmaName}

\newcounter {remark}

\newcounter {example}

\newcounter {definition}

\newcounter {problem}

\newcommand{\proofName}{Proof}
\newcommand{\answerName}{Answer}
\newcommand{\hintName}{Hint}

\newenvironment {proof} {\par\medskip\noindent{\bf \proofName\,\,}}
{\qed\par}

\def \Real {{I\!\!R}}
\def \Integer {{Z\!\!\!Z}}

\def \lnorm#1\rnorm {\vphantom{#1}\left\|\smash{#1}\right\|}
\def \lmod#1\rmod {\vphantom{#1}\left|\smash{#1}\right|}

\newcommand \bydef {\stackrel{\mbox{\scriptsize def}}{=}}
\newcommand \qed {\,\rule[-.23ex]{1.6ex}{1.6ex}}

\newcommand \name[1] {\mathop{\rm #1}\nolimits}

\renewcommand \phi {\varphi}
\renewcommand \rho {\varrho}

\author{Yu.M.Burman\thanks{Mathematical College of the Independent
University of Moscow}{\space}\thanks{e-mail: yurii@burman.mccme.ru}}
\title{Whitney's index formula in higher dimensions and Laplace integrals}
\date{March 18, 1997}
\begin{document}
\maketitle

\newcommand{\pder}[2] {\frac{\partial #1}{\partial #2}}
\newcommand{\pdertwo}[3] {\frac{\partial^2 #1}{\partial #2 \partial #3}}

\section{Introduction}\label{Sec:Intro}

Let $f = (f_1,f_2):\Real \to \Real^2$ be a {\em long curve}, i.e. a smooth
immersion such that $f_1(x) = x$ and $f_2(x) = 0$ for $x \le -1$ and $x \ge
1$. By the famous theorem of Whitney (see \cite{Index}) long curves are
classified by a single integral invariant, an index (rotation number).
Informally saying, index is a number of full rotations made by the tangent
vector $D_f(x) = (f_1'(x), f_2'(x)) \ne 0$ as $x$ goes from $-\infty$ to
$+\infty$ (or, equally, from $-1$ to $1$). The explicit formula for the
index is as follows:
\begin{equation}\label{Eq:Index1}
I(f) = \int_{-\infty}^\infty \frac{f_1''(x)f_2'(x) - f_2''(x)f_1'(x)}
{\bigl(f_1'(x)\bigr)^2 + \bigl(f_2'(x)\bigr)^2}\, dx
\end{equation}

The integral in the right-hand side of (\ref{Eq:Index1}) is $\int_\Real
(D_f)^* \omega$ where
\begin{equation}\label{Eq:Basis1}
\omega(z_1,z_2) = \frac{z_2 dz_1 - z_1 dz_2}{z_1^2 + z_2^2}
\end{equation}
is the closed form in $\Real^2 \setminus \{0\}$ forming a basis in its de
Rham cohomologies $H^1(\Real^2 \setminus \{0\}) = \Real$.

Suppose now that $f$ is generic. By Thom's transversality theorem (see
\cite{Singul}) $f$ has only a finite number of self-intersection points.
All these points are simple (i.e. only two branches of the curve meet) and
transversal. A point $a \in \Real^2$ of simple transversal
self-intersection can be equipped with a sign $\sigma(a) = \pm 1$ by the
following rule.  Let $f(x^{(1)}) = f(x^{(2)}) = a$ and $x^{(1)} < x^{(2)}$.
By transversality, the tangent vectors $D_f(x^{(1)})$ and $D_f(x^{(2)})$
are not parallel, and therefore form a basis in $\Real^2$. Choose an
orientation of $\Real^2$ and take
\begin{equation}\label{Eq:Sign1}
\sigma(a) = \left\{\begin{array}{ll}
+1, &\mbox{if the basis $\bigl(D_f(x^{(1)}),D_f(x^{(2)})\bigr)$}\\
& \mbox {gives a chosen orientation of $\Real^2$}, \\
-1, &\mbox{if the basis $\bigl(D_f(x^{(1)}),D_f(x^{(2)})\bigr)$} \\
&\mbox{gives an opposite orientation}.
\end{array}\right.
\end{equation}

Then there holds Whitney's index formula (see \cite{Index}):
\begin{equation}\label{Eq:SumSign1}
I(f) = \sum_a \sigma(a)
\end{equation}
where the sum is taken over all the self-intersection points of the curve
$f$.

The paper is devoted to the generalization of
(\ref{Eq:Index1})--(\ref{Eq:SumSign1}) to the case of smooth immersions
$f:\Real^n \to \Real^{2n}$ such that $f(x) = (x,0)$ for $\lmod x\rmod$
large. In Section \ref{Sec:SelfInt} we define an index and prove a ``sum of
signs'' formula for generic immersions. Section \ref{Sec:Integral} contains
a proof of an explicit integral formula for the index (for $n$ even); the
proof makes use of asymptotic expansions of certain Laplace integrals. The
formula obtained looks like $I(f) = \int_{\Real^n} (D_f)^* \omega$ where
$\omega$ is an $n$-form on the Stiefel variety $V(n,2n)$ and $D_f:\Real^n
\to V(n,2n)$ is a differential. It is proved in Section \ref{Sec:Stiefel}
that $\omega$ is a generator of the de Rham cohomology group
$H^n(V(n,2n))$.

\section{Index and self-intersection points}\label{Sec:SelfInt}

Recall that the map $f = (f_1, \dots, f_{2n}):\Real^n \to \Real^{2n}$ is
called an {\em immersion} if for any $x \in \Real^n$ the $(n \times
2n)$-matrix
\begin{equation}\label{Eq:DefDiff}
D_f(x) = \left(\pder{f_i}{x_j}\right)
\end{equation}
has the rank $n$ (i.e., belongs to the Stiefel variety $V(n,2n)$). We will
always suppose that the immersion $f$ is ``fixed at infinity'', i.e. $f(x)
= (x,0)$ as soon as $x$ lies outside the unit cube of $\Real^n$. Generic
immersions $\Real^n \to \Real^{2n}$ fixed at infinity have only a finite
number of self-intersection points, all of which are simple and
transversal.

The homotopy classes of immersions $\Real^n \to \Real^{2n}$ form an abelian
group $\name{Imm}_n$. Definition of the multiplication resembles that of a
homotopy group: take
\begin{displaymath}
\bigl(f_1*f_2\bigr)(x) = \left\{\begin{array}{ll}
f_1(2x_1 + 1, x_2, \dots, x_n), &\mbox{if $x_1 \le 0$}, \\
f_2(2x_1 - 1, x_2, \dots, x_n), &\mbox{if $x_1 \ge 0$}.
\end{array}\right.
\end{displaymath}

One of the formulations of the famous Smale classification theorem (see
\cite{ClassImm}) is
\begin{theorem}\label{Th:ClassImm}
The mapping $D_f:\name{Imm}_n \to \pi_n(V(n,2n))$ is a group isomorphism.
\end{theorem}

So to classify immersions we need to know the group $\pi_n(V(n,2n))$. The
answer is the following well-known lemma (see, for example,
\cite{FuchsFom}, page 173):

\begin{lemma}\label{Lm:HomStie}
\begin{eqnarray*}
\pi_i(V(n,2n)) &=& \begin{array}{ll}
0 &\mbox{for $1 \le i \le n-1$}.
\end{array}
\\
\pi_n(V(n,2n)) &=& \left\{\begin{array}{ll}
\Integer &\mbox{for $n=1$ and $n$ even},\\
\Integer_2 &\mbox{for odd $n \ge 3$}.
\end{array}\right.
\end{eqnarray*}
\end{lemma}

Thus, smooth immersions $f:\Real^n \to \Real^{2n}$ are classified by a
single invariant, an index $I(f)$ which is either an integer (for $n = 1$
and $n$ even) or an element of $\Integer_2$ (for odd $n \ge 3$). The case
$n = 1$ is covered by the results of Whitney, and in what follows we will
always assume that $n > 1$.

For $n$ even, a point $a \in \Real^{2n}$ of simple transversal
self-intersection of immersion $f$ can be equipped with a sign $\sigma(a) =
\pm 1$ by the following rule: if $f(x^{(1)}) = f(x^{(2)}) = a$ (and
$x^{(1)} \ne x^{(2)}$), then take
\begin{equation}\label{Eq:SignN}
\sigma(a) = \name{sgn} \det \left(\begin{array}{c}
D_f(x^{(1)})\\
D_f(x^{(2)})
\end{array}\right)
\end{equation}
(the right-hand side contains a determinant of the matrix $2n \times 2n$).
Geometrically this means that we choose positively oriented bases in the
tangent spaces to the image of $f$ at points $x^{(1)}$ and $x^{(2)}$. Since
$f$ is immersion, these $2n$ vectors form a basis in $\Real^{2n}$, and
$\sigma(a)$ is taken $\pm 1$ depending on the orientation of this basis
(orientation of $\Real^{2n}$ is fixed beforehand). Since $n$ is supposed to
be even, the result will not depend on the order in which the preimages
$x^{(1)}, x^{(2)} \in f^{-1}(a)$ are taken (for $n=1$ we supposed that
$x^{(1)} < x^{(2)}$, but there is no way to choose this order for $n > 1$).

\begin{theorem}\label{Th:SumSign}
\begin{equation}\label{Eq:SumSignN}
I(f) = \sum_a \sigma(a)
\end{equation}
where the sum is taken over all the self-intersection points of $f$, and
the summation is performed in $\Integer$ for $n$ even and in $\Integer_2$
for odd $n \ge 3$.
\end{theorem}

\begin{proof}
Let $f_t$ be a generic continuous family of smooth immersions $\Real^n \to
\Real^{2n}$ fixed at infinity. By Thom's transversality theorem the set of
self-intersection points of $f_t$ is changes continuously with $t$, except
for a finite number of points $t_1, \dots, t_k$ where two ``catastrophes''
may occur: either two self-intersection points merge and disappear, or,
conversely, emerge. It can be easily observed that for $n$ even the pair of
points to emerge or to perish always have opposite signs. It follows from
this that (for any $n$) the right-hand side of (\ref{Eq:SumSignN}) depends
only on the homotopic class of the immersion $f$. Thus, this right-hand
side defines a mapping $S$ from $\name{Imm}_n$ to $\Integer$ or
$\Integer_2$, depending on parity of $n$.

Definition of multiplication in $\name{Imm}_n$ shows that $S$ is a
homomorphism. By Theorem \ref{Th:ClassImm}, to prove equation
(\ref{Eq:SumSignN}) it is enough to show that $S$ is an isomorphism, i.e.
that its image is the whole $\Integer$ or $\Integer_2$. In other words it
should be proved that $1 \in \name{Im} S$, i.e. that there exists an
immersion $f^{(n)}: \Real^n \to \Real^{2n}$ with only one self-intersection
point.

Such immersion is easily constructed by induction on $n$. For $n = 1$ take
a smooth plane curve $f^{(1)}$ such that $f^{(1)}(-1/2) = f^{(1)}(1/2)$,
and there are no more self-intersection points. Suppose now that immersion
$f^{(n)}$ is already constructed, and the self-intersection is
$f^{(n)}(-1/2, 0, \dots, 0) = f^{(n)}(1/2, 0, \dots, 0)$. Let $\phi(t)$ be
a smooth function such that $\phi(-1/2) \ne \phi(1/2)$ and $\phi(t) = 0$
for $\lmod t\rmod \ge 1$. Then it can be easily seen that the
formula
\begin{eqnarray*}
f^{(n+1)}(x_1, \dots, x_{n+1}) &=& \bigl(f^{(n)}_1(x_1, \dots, x_n), \dots,
f^{(n)}_n(x_1, \dots, x_n), x_{n+1}, \\
&&f^{(n)}_{n+1}(x_1, \dots, x_n), \dots, f^{(n)}_{2n}(x_1, \dots, x_n),
x_{n+1} \phi(x_1)\bigr)
\end{eqnarray*}
defines an immersion $\Real^{n+1} \to \Real^{2n+2}$ whose sole
self-intersection is \linebreak $f^{(n+1)}(-1/2, 0, \dots, 0) =
f^{(n+1)}(1/2, 0, \dots, 0)$, so that induction step is made.
\end{proof}

\section{An integral formula for the index}\label{Sec:Integral}

Let us derive an explicit integral formula for the index $I(f)$. For the
reasons explained in the previous Section we restrict ourselves to the case
of $n$ even (an element of $\Integer_2$ hardly can appear as an integral).

Start with a technical result. Let $f:\Real^n \to \Real^{2n}$ be an
arbitrary smooth map (even not immersion) fixed at infinity. Consider a
mapping $A_f:\Real^{2n} \to \Real^{2n}$ given by the formula $A_f(x,y) =
f(x) - f(y)\ (x,y \in \Real^n)$ and recall that {\em degree of a smooth
mapping} $g:M \to N$ ($M$ and $N$ being oriented manifolds of the same
dimension) is defined as
\begin{displaymath}\label{Eq:DefDeg}
\name{deg} g = \sum_{g(x) = y} \name{sgn}(x).
\end{displaymath}
Here $y \in N$ is a generic point, and $\name{sgn}(x)$ is taken $1$ or $-1$
depending on whether the mapping $Df:T_x M \to T_y N$ preserves or inverts
orientation. The right-hand side of this equation can be called
multiplicity of $g$ at $y$.

\begin{lemma}\label{Lm:DegZero}
The degree of the mapping $A_f$ is zero.
\end{lemma}
\begin{proof}
Prove first that the degree of $A_f$ is defined. To do this denote $E
\subset \Real^{2n}$ the subspace of vectors $(a_1, \dots, a_n, 0, \dots,
0)$. Since $n > 1$, any two points $z_0, z_1 \in \Real^{2n} \setminus E$
can be connected with the path $z_t$ that does not intersect $E$. Notice
now that since $f$ is fixed at infinity, for every compact set $B \subset
\Real^{2n} \setminus E$ its full preimage $A_f^{-1}(B)$ is compact, and
therefore multiplicity of $A_f$ at $z_t$ does not depend on $t$. So,
multiplicity is the same for all generic points, and $\name{deg} A_f$ is
defined. But then $\name{deg} A_f = 0$ because the image of $A_f$ is not
the whole $\Real^{2n}$ (for example, if $z \in \Real^{2n}$ is large enough
and normal to the subspace $E$ then certainly $A_f^{-1}(z) = \emptyset$).
\end{proof}

Consider not $A_f$ as a change of variable in $\Real^{2n}$ and apply it to
the $2n$-form $\nu = \exp(-\lambda \lmod z\rmod^2/2) dz_1 \wedge \dots
\wedge dz_{2n}$. Using Lemma \ref{Lm:DegZero}, obtain the following result:
\begin{eqnarray}
J(f) &=& \int_{\Real^n \times \Real^n} \exp\left(-\lambda \lmod f(x) -
f(y)\rmod^2/2\right) \det
\left(\begin{array}{c}
D_f(x)\\
D_f(y)
\end{array}\right)\, dxdy \nonumber\\
&=& \int_{\Real^{2n}} A_f^* \nu = \name{deg} A_f \int_{\Real^{2n}} \nu = 0.
\label{Eq:IntIsZero}
\end{eqnarray}

Our idea is now to find the asymptotics of the integral
(\ref{Eq:IntIsZero}) for $\lambda \to +\infty$ using the general results
about asymptotics of Laplace integrals. Then we consider the principal term
of the asymptotics and equalize it to zero.

To formulate the result introduce first some notations. Let $\phi =
(\phi^j_i) \in V(n,2n)$ (i.e. $\name{Rk}\phi = n,\, 1 \le j \le 2n, 1 \le
i \le n$). For an arbitrary set of integers $J = \{j_1, \dots, j_n\}$ such
that $1 \le j_1 < \dots < j_n \le 2n$ denote
\begin{eqnarray*}
&&\mu(J) = n(n-1)/2 + j_1 + \dots + j_n, \\
&&{\frak M}_J \bydef \det \bigl(\phi^{\alpha_s}_s\bigr)_{\small
\begin{array}{l} s = 1, \dots, n \\ \alpha_s \notin J \end{array}}
\end{eqnarray*}
Denote also
\begin{eqnarray*}
&&U \bydef \left(\sum_{k=1}^{2n} \phi^k_i\phi^k_j\right),\\
&&(u_{ij}(\phi)) \bydef U^{-1}, \\
&&\det U = 1/u^2(\phi).
\end{eqnarray*}
(by the general theorem about Gram determinant, $\det U$ equals the sum of
squares of the $n$-th order minors of the matrix $\phi$, and is therefore
positive).

\begin{theorem}\label{Th:IntIndex}
Let $n$ be even, $f:\Real^n \to \Real^{2n}$ be a smooth immersion fixed at
infinity, and $D_f: \Real^n \to V(n,2n)$ be its differential. Then the
index of $f$ is given by the formula
\begin{equation}\label{Eq:IndexForm}
I(f) = \int_{\Real^n} (D_f)^* \omega
\end{equation}
where the $2n$-form $\omega$ at the point $\phi = (\phi_i^j) \in V(n,2n)$
is
\begin{eqnarray}
\omega(\phi) &=& - \frac{1}{2^{n+1}\pi^{n/2}(n/2)!} u(\phi) \sum_{i_1,
\dots, i_n = 1}^n u_{i_1i_2}(\phi) \dots u_{i_{n-1}i_n}(\phi) \nonumber\\
&\times& \sum_{\sigma \in \Sigma_n} \hskip-0.5cm \sum_{\small
\begin{array}{c} J = \{j_1, \dots, j_n\} \\ 1 \le j_1 < \dots < j_n \le 2n
\end{array}} \hskip-1cm (-1)^{\mu(J)} {\frak
M}_J(\phi)\, d\phi^{j_1}_{i_{\sigma(1)}} \wedge \dots \wedge
d\phi^{j_n}_{i_{\sigma(n)}}.\label{Eq:FormStie}
\end{eqnarray}
($\Sigma_n$ means the symmetric group).
\end{theorem}

To prove Theorem \ref{Th:IntIndex} we will need the following asymptotical
expansion found in \cite{Asymp} (Statement 4.1 at page 125):

\begin{lemma}\label{Lm:Expand}
Let $a,S$ be smooth real-valued functions at $\Real^k$, $S$ achieving its
minimum at the origin (and only at the origin). Let the Hessian $G =
\pdertwo{S}{u_i}{u_j}(0)$ be nondegenerate.  Denote $L$ the 2nd order
differential operator $L = \sum_{i,j=1}^k \bigl(G^{-1}\bigr)_{ij}
\pdertwo{}{u_i}{u_j}$, and $R(u) \bydef S(u) - S(0) - \frac12
\bigl(Gu,u\bigr)$. Then for $\lambda \to +\infty$ there is an aymptotic
expansion up to the $O(\lambda^{-\infty})$:
\begin{eqnarray}
\int_{\Real^k} a(u)\exp(-\lambda S(u))\, du \sim \exp(-\lambda S(0))
\left(\frac{2\pi}{\lambda}\right)^{k/2} \lmod \det G\rmod^{-1/2}
\nonumber\\
\times \sum_{p=0}^\infty \frac{1}{p!(2\lambda)^p} \left.L^p
\biggl(a(u)\exp(-\lambda R(u))\biggr)\right|_{u = 0}.\label{Eq:Expand}
\end{eqnarray}
\end{lemma}

{\def \proofName {Proof of Theorem \ref{Th:IntIndex}}
\begin{proof}
The left-hand side $J(f)$ of (\ref{Eq:IntIsZero}) is exactly the Laplace
integral mentioned in Lemma \ref{Lm:Expand} with $k = 2n$, $u = (x,y)$,
$S(x,y) = \lmod f(x) - f(y)\rmod^2/2$, and $a(x,y) = \det
\left(\begin{array}{c} D_f(x)\\ D_f(y) \end{array}\right)$. Minimum point
of $S$ is not unique but the general localization principle for Laplace
integrals (see \cite{Asymp}) tells that one should consider expansions
(\ref{Eq:Expand}) for all such points and sum them up (points other then
local minima make $O(\lambda^{-\infty})$ contributions).

Local minima of $S(x,y)$ are:
\begin{enumerate}
\item Self-intersection points of $f$: $f(x) = f(y)$ and $x \ne y$.
\item Diagonal points: $x = y$.
\item Other (``sporadic'') minima.
\end{enumerate}

At sporadic minima $s = S(x,y) > 0$. So expansion (\ref{Eq:Expand})
contains the term $\exp(-\lambda s)$ and therefore is already
$O(\lambda^{-\infty})$. Thus, sporadic minima will be neglected.

To self-intersection points $(x,y)$ Lemma \ref{Lm:Expand} can be applied
directly. The principal term here is $p = 0$. Denote $g_i\ (i = 1, \dots,
2n)$ the vector \linebreak $(\pder{f_i}{x_1}(x), \dots,
\pder{f_i}{x_n}(x),\\ -\pder{f_i}{y_1}(y), \dots, -\pder{f_i}{y_n}(y))$ and
observe that $G_{ij} = \bigl(g_i,g_j\bigr)$. By the general theorem about
Gram determinant we have $\det G = \left|\det \left(\begin{array}{c}
D_f(x)\\ D_f(y) \end{array}\right)\right|^2$, and the principal term is
\begin{displaymath}
\left(\frac{2\pi}{\lambda}\right)^n \lmod \det G\rmod^{-1/2} a(x,y) =
\left(\frac{2\pi}{\lambda}\right)^n \name{sgn} \det \left(\begin{array}{c}
D_f(x)\\ D_f(y) \end{array}\right).
\end{displaymath}
It follows now from Theorem \ref{Th:SumSign} that total contribution to the
principal term of asymptotics from all the self-intersection points is
equal to
\begin{equation}\label{Eq:ContSInt}
{\cal C}_{\mbox{\scriptsize self-intersections}} =
2\left(\frac{2\pi}{\lambda}\right)^n I(f).
\end{equation}
The factor $2$ is present because every self-intersection point appears
twice as a minimum of $S$: the first time as $(x,y)$ and the second time as
$(y,x)$.

Consider now the contribution from the diagonal. The minimum of $S$ here is
not isolated (and therefore degenerate) but it is easy to see that we
should consider asymptotic expansion of the integral over $y$ (with $x$
fixed) and then integrate it over $x$ (necessary theorems about uniformity
of asymptotic expansion (\ref{Eq:Expand}) are contained in \cite{Asymp}).
If $x$ is fixed then $y = x$ is a nondegenerate minimum of $S$, and $G_{ij}
= \bigl(h_i,h_j\bigr)$ where $h_i = (\pder{f_1}{x_i}(x), \dots,
\pder{f_{2n}}{x_i}(x))$. By the general theorem about Gram determinant
\begin{eqnarray}
g(x) \bydef \lmod \det G\rmod^{-1/2} &=& \left|\det \left(\sum_{k=1}^{2n}
\pder{f_k}{x_i}(x) \pder{f_k}{x_j}(x)\right)\right|^{-1/2} \nonumber\\
&=& \left(\sum_{1 \le i_1 < \dots < i_n \le 2n} \left|\det
\left(\pder{f_{i_s}}{x_s}(x)\right)\right|^2\right)^{-1/2}.\label{Eq:Def_g}
\end{eqnarray}
The right-hand side contains the sum of squares of all the $n$-th order
minors of the $(2n \times n)$-matrix $\left(\pder{f_i}{x_j}(x)\right)$.

First we are to find the leading term of expansion. Let $v(t), w(t)$ be
smooth functions of one variable. Then it is easy to see that
\begin{eqnarray}
\frac{d^m}{dt^m} &&\biggl(v(t)\exp (\lambda w(t))\biggr) = \exp (\lambda
w(t))\sum_u \lambda^u \nonumber\\
&\times& \sum_{\small \begin{array}{c}
r, k_1, \dots, k_u \ge 0, \\
r+k_1+\dots+k_u = m
\end{array}}
b_{r,k_1, \dots, k_u} v^{(r)}(t) w^{(k_1)}(t)\dots
w^{(k_u)}(t)\label{Eq:DiffExp}
\end{eqnarray}
for some integers $b_{r,k_1, \dots, k_u} \ge 0$ (which we need not
specify). Apply this formula with $v = a,\ w = R$ to expansion
(\ref{Eq:Expand}). Here $a(x,y) = O(\lmod x-y\rmod^n)$ and $R(x,y) =
O(\lmod x-y\rmod^3)$, and the degree of the operator $L$ is $2$ (so, $L^p$
contains differentiations of the order $2p$). Thus, the $p$-th term of
expansion will contain $\lambda^{u-p-n/2}$ where
\begin{eqnarray*}
&&2p = r + k_1 + \dots + k_u, \\
&&r \ge n, \space k_i \ge 3.
\end{eqnarray*}
Thus $2p \ge n + 3u$, and therefore $u-p \le -n/2$. So, the leading term is
\begin{equation}\label{Eq:LeadTerm}
\frac{\pi^{n/2}}{(n/2)!} \frac{1}{\lambda^n}\int_{\Real^n} g(x)
\left.(L^{n/2} a(x,y))\right|_{y=x}\, dx
\end{equation}
where $g(x)$ is given by (\ref{Eq:Def_g}).

Denote $G^{-1} \bydef \bigl(g_{ij}(x)\bigr)$. Then, by definition,
\begin{displaymath}
L = \sum_{i,j = 1}^n g_{ij}(y) \pdertwo{}{y_i}{y_j}
\end{displaymath}
and
\begin{displaymath}
L^{n/2} = \sum_{i_1, \dots, i_n = 1}^n g_{i_1i_2}(y) \dots
g_{i_{n-1}i_n}(y) \frac{\partial^n}{\partial y_{i_1} \dots \partial
y_{i_n}}.
\end{displaymath}
This yields:
\begin{eqnarray*}
\left.(L^{n/2} a(x,y))\right|_{y=x} &=& \sum_{i_1, \dots, i_n = 1}^n
g_{i_1i_2}(x) \dots g_{i_{n-1}i_n}(x) \left.\frac{\partial^n}{\partial y_{i_1}
\dots \partial y_{i_n}} \det \left(\begin{array}{c} \pder{f_j}{x_k}(x) \\
\pder{f_j}{y_k}(y)\end{array}\right)\right|_{y=x} \\
&=& \sum_{i_1, \dots, i_n = 1}^n g_{i_1i_2}(x) \dots g_{i_{n-1}i_n}(x)
\sum_{\sigma \in \Sigma_n} \det \left(\begin{array}{c} \pder{f_j}{x_k}(x) \\
\pdertwo{f_j}{x_k}{x_{i_{\sigma(k)}}}(x)\end{array}\right).
\end{eqnarray*}
Here $\Sigma_n$ means the symmetric group. The subscript $j$ in
determinants runs from $1$ to $2n$, the subscript $k$, from $1$ to $n$.
Thus, total contribution to the principal term of the asymptotics of $J(f)$
from the diagonal $y = x$ equals:
\begin{equation}\label{Eq:ContrDiag}
\begin{array}{l}
{\cal C}_{\mbox{\scriptsize diagonal}} = \frac{\pi^{n/2}}{(n/2)!}
\frac{1}{\lambda^n}\int_{\Real^n} g(x) \sum_{i_1, \dots, i_n = 1}^n
g_{i_1i_2}(x) \dots g_{i_{n-1}i_n}(x) \times \\
\hphantom{{\cal C}_{\mbox{\scriptsize diagonal}} =
\frac{\pi^{n/2}}{(n/2)!}} \times \sum_{\sigma \in \Sigma_n} \det
\left(\begin{array}{c} \pder{f_j}{x_k}(x) \\
\pdertwo{f_j}{x_k}{x_{i_{\sigma(k)}}}(x)\end{array}\right).
\end{array}
\end{equation}

But it follows from (\ref{Eq:IntIsZero}) that ${\cal C}_{\mbox{\scriptsize
diagonal}} + {\cal C}_{\mbox{\scriptsize self-intersections}} = 0$. Thus,
combining (\ref{Eq:ContSInt}) and (\ref{Eq:ContrDiag}), one has
\begin{equation}\label{Eq:IndThrInt}
\begin{array}{l}
I(f) = -\frac{1}{2^{n+1}\pi^{n/2}(n/2)!} \int_{\Real^n} g(x) \sum_{i_1,
\dots, i_n = 1}^n g_{i_1i_2}(x) \dots g_{i_{n-1}i_n}(x) \times \\
\hphantom{I(f) = -\frac{1}{2^{n+1}\pi^{n/2}(n/2)!}} \times \sum_{\sigma \in
\Sigma_n} \det \left(\begin{array}{c} \pder{f_j}{x_k}(x) \\
\pdertwo{f_j}{x_k}{x_{i_{\sigma(k)}}}(x)\end{array}\right).
\end{array}
\end{equation}

Notice now that for any $1 \le j_1, \dots, j_n \le 2n,\ 1 \le i_1, \dots,
i_n \le n$
\begin{equation}\label{Eq:SubstDet}
\det \left(\begin{array}{c} \pdertwo{f_{j_s}}{x_s}{x_{i_s}}
\end{array}\right)\, dx_1 \wedge \dots \wedge dx_n = (D_f)^*
d\phi^{j_1}_{i_1} \wedge \dots \wedge d\phi^{j_n}_{i_n}
\end{equation}
Now to get (\ref{Eq:IndexForm}) expand the determinant in
(\ref{Eq:IndThrInt}) in the first $n$ rows.
\end{proof}
}

\section{An explicit formula for the generator of
$H^n(V(n,2n))$}\label{Sec:Stiefel}

Take up again Lemma \ref{Lm:HomStie}. By Gurevich theorem (see e.g.
\cite{FuchsFom}) for $n$ even, $H^n_{DR}(V(n,2n)) = \Real$ is the first
nontrivial cohomology group of the Stiefel variety $V(n,2n)$. Let
$\omega_0$ be a closed $n$-form at $V(n,2n)$ whose integral over the
fundamental cycle (generator of the $H_n(V(n,2n),\Integer)$) is $1$. Then
the homotopic class $[\gamma] \in \pi_n(V(n,2n)) = \Integer$ of the spheroid
$\gamma:S^n \to V(n,2n)$ equals to $\int_{S^n} \gamma^* \omega_0$. Now
Theorems \ref{Th:IntIndex} and \ref{Th:ClassImm} suggest that $\omega =
\omega_0$, i.e.,

\begin{theorem}\label{Th:OmeBasis}
The cohomology class of the $n$-form $\omega$ given by (\ref{Eq:FormStie})
forms the basis in $H^n_{DR}(V(n,2n))$. The value of this class on the
generator of the $H_n(V(n,2n),\Integer)$ is $1$.
\end{theorem}

The main step in the proof of Theorem \ref{Th:OmeBasis} is the following
\begin{lemma}\label{Lm:OmeClos}
The form $\omega$ is closed.
\end{lemma}

\begin{proof}
By Theorem \ref{Th:IntIndex} the equality
\begin{equation}\label{Eq:IntOmOm0}
\int_\Real^n \gamma^*\omega = \int_\Real^n \gamma^*\omega_0
\end{equation}
holds for any {\em holonomic} $\gamma = (\gamma^i_j): \Real^n \to V(n,2n)$.
The holonomity means that $\gamma = D_f$ for some immersion $f:  \Real^n
\to \Real^{2n}$ fixed at infinity, or, equivalently, that the equations
\begin{equation}\label{Eq:HolGamma}
\pder{\gamma^i_j(x)}{x_k} = \pder{\gamma^i_k(x)}{x_j} \quad \mbox{for any
$i,j,k$}
\end{equation}
and
\begin{equation}\label{Eq:InfGamma}
\gamma^i_j(x) = \left\{\begin{array}{rl}
1, &i = j \le n, \\
0, &\mbox{otherwise.}
\end{array}\right.
\end{equation}
(for $x$ lying outside the unit cube of $\Real^n$) hold.

Consider a one-parameter family of mappings $\gamma_t = (\gamma^i_j +
t\delta^i_j)$. The definition of Lie derivative $\cal L$ and the Cartan
formula ${\cal L}_X = \iota_X d + d \iota_X$ imply that
\begin{eqnarray}
\frac{d}{dt} \int_{\Real^n} \gamma_t^* \omega &=& \sum_{i=1}^{2n}
\sum_{j=1}^n \int_{\Real^n} \delta^i_j \gamma_t^* {\cal
L}_{\partial/\partial \phi^i_j} \omega \nonumber\\
&=& \sum_{i=1}^{2n} \sum_{j=1}^n \int_{\Real^n} \delta^i_j \gamma_t^*
\iota_{\partial/\partial \phi^i_j} d\omega.\label{Eq:Deriv}
\end{eqnarray}

Taking in (\ref{Eq:Deriv}) $t = 0$ and $\gamma = \gamma_0 = D_f$, and
taking (\ref{Eq:IntOmOm0})--(\ref{Eq:InfGamma}) into consideration, one
obtains the equality
\begin{equation}\label{Eq:SumDelta}
\sum_{i=1}^{2n} \sum_{j=1}^n \int_{\Real^n} \delta^i_j (D_f)^*
\iota_{\partial/\partial \phi^i_j} d\omega = 0
\end{equation}
for every $f$ and every $\delta = (\delta^i_j(x))$ such that
\begin{equation}\label{Eq:DeltHol}
\pder{\delta^i_j(x)}{x_k} = \pder{\delta^i_k(x)}{x_j} \quad \mbox{for any
$i,j,k$}
\end{equation}
and
\begin{equation}\label{Eq:DeltInf}
\delta^i_j(x) = 0
\end{equation}
for $x$ lying outside the unit cube of $\Real^n$.

Take some immersion $f: \Real^n \to \Real^{2n}$ fixed at infinity, and
denote $P^i_j(x)$ the function such that $P^i_j(x)\, dx_1 \wedge \dots
\wedge dx_n = (D_f)^* \iota_{\partial/\partial \phi^i_j} d\omega$. Denote
also $Q^i_{jk}$ a function such that $P^i_j = \pder{Q^i_{jk}}{x_k}$.
$Q^i_{jk}$ is defined uniquely up to addition of a function $C^i_{jk}$
independent of $x_k$. If (\ref{Eq:DeltHol}) and (\ref{Eq:DeltInf}) are
fulfilled then
\begin{eqnarray}
\sum_{i=1}^{2n} \sum_{j=1}^n \int_{\Real^n} \delta^i_j P^i_j\, dx &=&
\frac{1}{n} \sum_{i=1}^{2n} \sum_{j,k=1}^n \int_{\Real^n} \delta^i_j
\pder{Q^i_j}{x_k}\, dx \nonumber\\
&=& - \frac{1}{n} \sum_{i=1}^{2n} \sum_{j,k=1}^n \int_{\Real^n} Q^i_j
\pder{\delta^i_j}{x_k}\,dx. \label{Eq:FunctDelt}
\end{eqnarray}

Since (\ref{Eq:FunctDelt}) is zero for $\delta$ satisfying
(\ref{Eq:DeltHol}) and (\ref{Eq:DeltInf}), the standard Riesz lemma implies
that
\begin{equation}\label{Eq:DiffPjk}
\pder{P^i_j}{x_k} + \pder{P^i_k}{x_j} = 0 \quad \mbox{for any $i,j,k$}.
\end{equation}
In particular, it holds for $j=k$, i.e.
\begin{displaymath}
\pder{P^i_j}{x_j} = 0,
\end{displaymath}
or
\begin{equation}\label{Eq:DiffDfj}
{\cal L}_{\partial/\partial x_j} (D_f)^* \iota_{\partial/\partial \phi^i_j}
d\omega = 0.
\end{equation}
Equation (\ref{Eq:DiffDfj}) and boundary condition (\ref{Eq:InfGamma})
imply that $(D_f)^* \iota_{\partial/\partial \phi^i_j} d\omega = 0$. Since
$f$ is arbitrary, it is possible only if $\iota_{\partial/\partial
\phi^i_j} d\omega = 0$ for any $i,j$, which means that $d\omega = 0$. Lemma
is proved.
\end{proof}

{\def \proofName {Proof of Theorem \ref{Th:OmeBasis}}
\begin{proof}
Since $\omega$ is closed and the cohomologies $H^n(V(n,2n))$ are
one-dimensional, there exist $\lambda \in \Real$ and the $(n-1)$-form $\nu$
such that $\omega = \lambda \omega_0 + d\nu$. Then for any $\gamma:\Real^n
\to V(n,2n)$ fixed at infinity one obtains
\begin{equation}\label{Eq:IntOmLOm0}
\int_\Real^n \gamma^*\omega = \lambda \int_\Real^n \gamma^*\omega_0
\end{equation}
Comparing (\ref{Eq:IntOmLOm0}) with (\ref{Eq:IntOmOm0}) one obtains
$\lambda = 1$. Theorem is proved.
\end{proof}
}

\end{document}